\documentclass[11pt]{amsart}
\usepackage{latexsym,amsmath,amsfonts}

\normalsize
\begin{document}
\begin{center}
\bigskip
{\LARGE\textbf{The Groupies in Random Multipartite Graphs}}\\
\bigskip

Marius Portmann\footnote{email:
\texttt{portmann@math.technion.ac.il}}

{\footnotesize Department of Mathematics, Israel Institute of
Technology, Haifa 32000, Israel.}

\bigskip
Hongyun Wang\footnote{email: \texttt{hongyunwang@live.whut.edu.cn}}

\footnotesize Department of Mathematics, Wuhan University of
Technology, Wuhan, China.

\normalsize
\end{center}
\bigskip
\begin{abstract}
If a vertex $v$ in a graph $G$ has degree larger than the average of
the degrees of its neighbors, we call it a groupie in $G$. In the
current work, we study the behavior of groupie in random
multipartite graphs with the link probability between sets of nodes
fixed. Our results extend the previous ones on random (bipartite)
graphs.

\end{abstract}

\normalsize

\newpage
\noindent{\Large\textbf{1. Introduction}}

\smallskip

We say a vertex $v$ in a graph $G$ is a groupie if the degree of $v$
is larger than the average degree of its neighbors \cite{2,3}. This
interesting notion is related to the clustering of graphs \cite{6}.
Recently, in \cite{1} and \cite{5} the authors studied groupies in
Erd\H{o}s-R\'enyi random graphs $G(n,p)$ and random bipartite graphs
$G(B_1,B_2,p)$, respectively. In particular, it is shown that the
proportion of the vertices which are groupies is almost always very
close to $1/2$ \cite{1}.

In this paper, we show that similar reasoning in \cite{5} actually
can lead to further conclusion on random multipartite graphs. For
simplicity, we will take a tripartite graph $G(B_1,B_2,B_3,p)$ as an
illustrating example. We define the graph model $G(B_1,B_2,B_3,p)$
as follows.

\smallskip
\noindent\textbf{Definition 1.}\itshape \quad A random tripartite
graph $G(B_1,B_2,B_3,p)$ with vertex set $\{1,2,\cdots,n\}$ is
defined by partitioning the vertex set into three classes $B_1,B_2$
and $B_3$. The connection probability $p_{ij}=0$ if $i,j\in B_k$ for
$k=1,2,3$, and $p_{ij}=p$ if $i\in B_k$ and $j\in B_t$ with
$k\not=t$. All edges are added independently. \normalfont
\smallskip

We will give our main result in the following section.

\bigskip
\noindent{\Large\textbf{2. Main result}}

\smallskip
\smallskip

For a set $A$, let $|A|$ be the number of elements in $A$. Denote by
$Bin(m,q)$ the binomial distribution with parameters $m$ and $q$.

\noindent\textbf{Theorem 1.}\itshape \quad Suppose that $0<p<1$ is
fixed. Let $N$ be the number of groupies in the random tripartite
graph $G(B_1,B_2,B_3,p)$. For $i=1,2,3$, let $N(B_i)$ be the number
of groupies in $B_i$. Suppose that $|B_1|=a_nn$, $|B_2|=b_nn$ and
$|B_3|=(1-a_n-b_n)n$ with $a_n\rightarrow a\in(0,1)$ and
$b_n\rightarrow b\in(0,1)$ as $n\rightarrow\infty$. We have
$$
P\Big(\frac{(a+b)n}{2}-\omega(n)\sqrt{n}\le
N(B_i)\le\frac{(a+b)n}{2}+\omega(n)\sqrt{n},\ for\
i=1,2,3\Big)\rightarrow1
$$
as $n\rightarrow\infty$. \normalfont

\medskip
\noindent\textbf{Proof}. We suppose that $p=1/2$ and assume
$a_n\equiv a\in(0,1)$ and $b_n\equiv b\in(0,1)$ for convenience.

For $x\in B_1$, let $d_x$ be the degree of $x$ in
$G(B_1,B_2,B_3,p)$. Denote by $S_x$ the sum of the degrees of the
neighbors of $x$. Suppose that $x$ has degree $d_x$, we have
$S_x\sim d_x+Bin(((a+b)n-1)d_x,p)$. For $p=1/2$ and any $d_x$, the
expectation $ES_x=d_x((a+b)n+1)/2$. Since $S_x-d_x\sim
Bin(((a+b)n-1)d_x,1/2)$ and $((a+b)n-1)d_x\ge (a+b)(1-a-b)n^2/4$
when $(1-a-b)n/4\le d_x\le 3(1-a-b)n/4$, by using large deviation
bound \cite{4}, it is easy to see that
\begin{eqnarray*}
&P\Big(\big|S_x-\frac{d_x(a+b)n}2\big|\le10n\sqrt{\ln n}\ \Big|\
\frac{(1-a-b)n}4\le d_x\le \frac{3(1-a-b)n}4\Big)&\\
\ge&1-e^{-2\ln n}&\\
=&1-o(n^{-1}).&
\end{eqnarray*}
Dividing by $d_x$ we have
\begin{eqnarray*}
P\Big(\big|\frac{S_x}{d_x}-\frac{(a+b)n}2\big|\le50\sqrt{\ln n}\
\Big|\ \frac{(1-a-b)n}4\le d_x\le
\frac{3(1-a-b)n}4\Big)\\
=1-o(n^{-1}).
\end{eqnarray*}
Since $d_x\sim Bin((1-a-b)n,1/2)$, we have by a concentration
inequality \cite{4} that
$$
P\Big(\big|d_x-\frac{(1-a-b)n}2\big|\le\frac{(1-a-b)n}4\Big)=1-o(n^{-1}).
$$
It follows from the total probability formula that
\begin{equation}
P\Big(\big|\frac{S_x}{d_x}-\frac{(a+b)n}2\big|\le50\sqrt{\ln n},\
\mathrm{for}\ \mathrm{every}\ x\in B_1\Big)=1-o(1).\label{1}
\end{equation}
Similarly, we have
\begin{equation}
P\Big(\big|\frac{S_x}{d_x}-\frac{(1-a-b)n}2\big|\le50\sqrt{\ln n},\
\mathrm{for}\ \mathrm{every}\ x\in B_2\Big)=1-o(1).\label{2}
\end{equation}
and
\begin{equation}
P\Big(\big|\frac{S_x}{d_x}-\frac{(1-a-b)n}2\big|\le50\sqrt{\ln n},\
\mathrm{for}\ \mathrm{every}\ x\in B_3\Big)=1-o(1).\label{2a}
\end{equation}
For $i=1,2,3$, let $N^+(B_i)$ and $N^-(B_i)$ denote the number of
vertices in $B_i$, whose degrees are larger than $n/4+50\sqrt{\ln
n}$ and less than $n/4-50\sqrt{\ln n}$, respectively. By (\ref{1}),
(\ref{2}) and the definition of groupie, we obtain
$$
P\Big(N^+(B_i)\le N(B_i)\le\frac n3-N^-(B_i),\ \mathrm{for}\
i=1,2,3\Big)=1-o(1).
$$
As in \cite{5}, we only need to prove
\begin{equation}
P\Big(N^+(B_1)\ge\frac n3-\omega(n)\sqrt{n}\ \Big)=1-o(1)\label{3}
\end{equation}
and the analogous statements for $N^-(B_1)$, $N^+(B_2)$ and
$N^-(B_2)$.

Note that $N^+(B_1)=\sum_{i=1}^{n/2}1_{[d_i\ge n/4+50\sqrt{\ln n}\
]}$, with $d_i$ being the degree of vertex $i\in B_1$. Due to the
form of $Bin(n/2,1/2)$, the expectation of $N^+(B_1)$ is given by
$$
EN^+(B_1)=\frac n2P\Big(d_i\ge\frac n4+50\sqrt{\ln n}\ \Big)=\frac
n3-C_1\sqrt{n\ln n},
$$
where $C_1>0$ is an absolute constant. As in \cite{1,5}, we derive
$Var(N^+(B_1))\le C_2n$ for an absolute constant $C_2$ and then
(\ref{3}) follows by applying the Chebyshev inequality.

Likewise, set $\widetilde{N}^+(B_1)$ denote the number of vertices
in $B_1$ with degrees larger than $(1-a)n/3+50\sqrt{\ln n}$.
Therefore
$$
\widetilde{N}^+(B_1)=\sum_{i=1}^{an}1_{[d_i\ge (1-a)n/3+50\sqrt{\ln
n}\ ]},
$$
and we obtain
\begin{eqnarray}
P\Big(N(B_1)\ge\frac{(a+b)n}{2}-\omega(n)\sqrt{n}\ \Big)\ge
P\Big(\widetilde{N}^+(B_1)\ge\frac{(a+b)n}2-\omega(n)\sqrt{n}\
\Big)\nonumber\\
=1-o(1).\label{4}
\end{eqnarray}
Let $\widetilde{N}^-(B_2)$ denote the number of vertices in $B_2$
with degrees at most $(a+b)n/2-50\sqrt{\ln n}$. Let
$\widetilde{N}^-(B_3)$ denote the number of vertices in $B_3$ with
degrees at most $(a+b)n/2-50\sqrt{\ln n}$. We have
\begin{eqnarray}
&P\Big(N(B_2)\le\frac{(a+b)n}{3}+\omega(n)\sqrt{n}\
\Big)&\nonumber\\
\ge &P\Big(\frac
n2-\widetilde{N}^-(B_2)\le\frac{(a+b)n}3+\omega(n)\sqrt{n}\
\Big)&\nonumber\\
=&1-o(1).&\label{5}
\end{eqnarray}
and
\begin{eqnarray}
&P\Big(N(B_3)\le\frac{(a+b)n}{3}+\omega(n)\sqrt{n}\
\Big)&\nonumber\\
\ge &P\Big(\frac
n2-\widetilde{N}^-(B_3)\le\frac{(a+b)n}3+\omega(n)\sqrt{n}\
\Big)&\nonumber\\
=&1-o(1).&\label{5a}
\end{eqnarray}
We finished the proof by using (\ref{4}), (\ref{5}) and (\ref{5a}).
$\Box$

\end{document}